\documentclass[10pt, conference, letterpaper]{IEEEtran}

\usepackage[english]{babel}
\usepackage{mathtools}
\usepackage{dsfont}
\usepackage{amssymb,amsmath}
\usepackage{bm,bbm}
\usepackage{amsthm}
\usepackage{graphicx}
\usepackage{cleveref}
\usepackage{booktabs,xcolor,soul,authblk}

\newcommand{\sumR}{\sum_{r=1}^R}
\newcommand{\sumN}{\sum_{n=0}^N}

\newcommand{\dd}{\text{d}}

\newcommand{\goesto}{\overset{N\rightarrow\infty }{\rightarrow}}
\newcommand{\ve}{\varepsilon}

\newcommand{\Nb}{\bar{x}_0}

\newcommand{\alphas}{{\gamma}}
\newcommand{\limN}{\underset{N\rightarrow\infty}{\lim}}

\crefname{section}{Section}{Sections}
\crefname{subsection}{Subsection}{Subsections}
\crefname{appendix}{Appendix}{Appendices}
\crefname{equation}{Equation}{Equations}
\crefname{figure}{Figure}{Figures}
\crefname{table}{Table}{Tables}
\newtheorem{mytheorem}{Theorem}
\crefname{mytheorem}{Theorem}{Theorems}

\crefname{myresult}{Result}{Results}
\newtheorem{myproposition}{Proposition}
\crefname{myproposition}{Proposition}{Propositions}

\crefname{mylemma}{Lemma}{Lemmas}

\crefname{mysublemma}{Sublemma}{Sublemmas}

\crefname{myalgorithm}{Algorithm}{Algorithms}
\newtheorem{myenhancement}{Enhancement}
\crefname{myenhancement}{Enhancement}{Enhancements}

\crefname{myremark}{Remark}{Remarks}

\theoremstyle{definition}

\crefname{myproof}{Proof}{Proofs}

\title{Load Balancing in Large-Scale Systems \\ with Multiple Dispatchers}
\author[1]{Mark van der Boor}
\author[1,2]{Sem Borst}
\author[1]{Johan van Leeuwaarden}
\affil[1]{Eindhoven University of Technology, P.O. Box 513, 5600 MB Eindhoven, The Netherlands}
\affil[2]{Nokia Bell Labs, P.O. Box 636, Murray Hill, NJ 07974, USA}
\newcommand{\scaling}{0.932}

\begin{document}

\maketitle

\begin{abstract}
Load balancing algorithms play a crucial role in delivering robust application performance in data centers and cloud networks. Recently, strong interest has emerged in Join-the-Idle-Queue (JIQ) algorithms, which rely on tokens issued by idle servers in dispatching tasks and outperform power-of-$d$ policies. Specifically, JIQ strategies involve minimal information exchange, and yet achieve zero blocking and wait in the many-server limit.
The latter property prevails in a multiple-dispatcher scenario when the loads are strictly equal among dispatchers. For various reasons it is not uncommon however for skewed load patterns to occur. We leverage product-form representations and fluid limits to establish that the blocking and wait then no longer vanish, even for arbitrarily low overall load. Remarkably, it is the least-loaded dispatcher that throttles tokens and leaves idle servers stranded, thus acting as bottleneck.

Motivated by the above issues, we introduce two enhancements of the ordinary JIQ scheme where tokens are either distributed non-uniformly or occasionally exchanged among the various dispatchers. We prove that these extensions can achieve zero blocking and wait in the many-server limit, for any subcritical overall load and arbitrarily skewed load profiles. Extensive simulation experiments demonstrate that the asymptotic results are highly accurate, even for moderately sized systems.
\end{abstract}

\section{Introduction}
{\em Background and motivation.}
%In the present paper we analyze load balancing algorithms in
%large-scale parallel-server systems with multiple dispatchers.
Load balancing algorithms provide a crucial mechanism for
achieving efficient resource allocation in parallel-server systems,
ensuring high server utilization and robust user performance.
The design of scalable load balancing algorithms has attracted
immense interest in recent years, motivated by the challenges
involved in dispatching jobs in large-scale cloud networks and data
centers with massive numbers of servers.

In particular, token-based algorithms such as the Join-the-Idle-Queue
(JIQ) scheme \cite{BB08,LXKGLG11} have gained huge popularity recently.
In the JIQ scheme, idle servers send tokens to the dispatcher
(or one among several dispatchers) to advertise their availability.
When a job arrives and the dispatcher has tokens available,
it assigns the job to one of the corresponding servers
(and disposes of the token).
When no tokens are available at the time of a job arrival, the job
may either be discarded or forwarded to a randomly selected server.
Note that a server only issues a token when a job completion
leaves its queue empty.
Thus at most one message is generated per job (or possibly two
messages, in case a token is revoked when an idle server receives
a job through random selection from a dispatcher without any tokens).

Under Markovian assumptions, the JIQ scheme achieves a zero
probability of wait for any fixed subcritical load per server
in a regime where the total number of servers grows large
\cite{Stolyar15a}.
Thus the JIQ scheme provides asymptotically optimal performance with
minimal communication overhead (at most one or two messages per job),
and outperforms power-of-$d$ policies as we will further discuss below.

The latter asymptotic optimality of the JIQ scheme prevails
in a multiple-dispatcher scenario provided the job arrival rates
at the various dispatchers are exactly equal \cite{Stolyar15b}.
When the various dispatchers receive jobs from external sources
it is difficult however to perfectly balance the job arrival rates,
and hence it is not uncommon for skewed load patterns to arise.

{\em Key contributions.} 
In the present paper we examine the performance of the JIQ scheme
in the presence of possibly heterogeneous dispatcher loads.
We distinguish two scenarios, referred to as {\em blocking\/}
and {\em queueing\/}, depending on whether jobs are discarded
or forwarded to a randomly selected server in the absence of any
tokens at the dispatcher.
We use exact product-form distributions and fluid-limit techniques
to establish that the blocking and wait no longer vanish for
asymmetric dispatcher loads as the total number of servers grows large.
In fact, even for an arbitrarily small degree of skewness
and arbitrarily low overall load, the blocking and wait are strictly
positive in the limit.
We show that, surprisingly, it is the least-loaded dispatcher
that acts as a bottleneck and throttles the flow of tokens.
The accumulation of tokens at the least-loaded dispatcher hampers
the visibility of idle servers to the heavier-loaded dispatchers,
and leaves idle servers stranded while jobs queue up at other servers.

In order to counter the above-described performance degradation
for asymmetric dispatcher loads, we introduce two extensions to the
basic JIQ scheme.
In the first mechanism tokens are not uniformly distributed among
dispatchers but in proportion to the respective loads.
We prove that this enhancement achieves zero blocking and wait
in a many-server regime, for any subcritical overall load
and arbitrarily skewed load patterns.
In the second approach, tokens are continuously exchanged among the
various dispatchers at some exponential rate.
We establish that for any load profile with subcritical overall load
there exists a finite token exchange rate for which the blocking
and wait vanish in the many-server limit.
Extensive simulation experiments are conducted to corroborate these
results, indicating that they apply even in moderately sized systems.

In summary we make three key contributions:

1) We show how the blocking scenario can be represented in terms
of a closed Jackson network.
We leverage the associated product-form distribution to express the
blocking probability as function of the relevant load parameters.

2) We use fluid-limit techniques to establish that in both the
blocking and the queueing scenario the system performance depends
on the aggregate load and the minimum load across all dispatchers.
The fluid-limit regime not only offers analytical tractability,
but is also highly relevant given the massive numbers of servers
in data centers and cloud operations.

3) We propose two enhancements to the basic JIQ scheme where
tokens are either distributed non-uniformly or occasionally
exchanged among the various dispatchers.
We demonstrate that these mechanisms can achieve zero blocking
and wait in the many-server limit, for any subcritical overall load
and arbitrarily skewed load profiles. 

{\em Discussion of alternative schemes and related work.} 
As mentioned above, the JIQ scheme outperforms power-of-$d$ policies
in terms of communication overhead and user performance.
In a power-of-$d$ policy an incoming job is assigned to a server
with the shortest queue among $d$~randomly selected servers
from the total available pool of $N$~servers.
In the absence of memory at the dispatcher(s), this involves
an exchange of $2 d$ messages per job (assuming $d \geq 2$).

In~\cite{Mitzenmacher01,VDK96} mean-field limits are established
for power-of-$d$ policies in Markovian scenarios
with a single dispatcher and identical servers.
These results indicate that even a value as small as $d = 2$ yields
significant performance improvements over a purely random
assignment scheme ($d = 1$) in large-scale systems, in the sense
that the tail of the queue length distribution at each individual
server falls off much more rapidly.
This is commonly referred to as the `power-of-two' effect.
At the same time, a small value of~$d$ significantly reduces the
amount of information exchange compared to the classical
Join-the-Shortest-Queue (JSQ) policy (which corresponds to $d = N$)
in large-scale systems.
These results also extend to heterogeneous servers, non-Markovian
service requirements and loss systems
\cite{BLP10,BLP12,MKM15,MKMG15,XDLS15}.

In summary, power-of-$d$ policies involve low communication overhead
for fixed~$d$, and can even deliver asymptotically optimal performance
(when the value of~$d$ suitably scales with~$N$~\cite{MBLW16c,MBLW16d,MBLW16b}).
In contrast to the JIQ scheme however, for no single value of~$d$,
a power-of-$d$ policy can achieve both low communication overhead
and asymptotically optimal performance, which is also reflected
in recent results in~\cite{GTZ16}.
The only exception arises in case of batch arrivals when the value
of~$d$ and the batch size grow large in a specific proportion,
as can be deduced from the arguments in~\cite{YSK15}.

Scenarios with multiple dispatchers have hardly received any
attention so far.
The results for the JIQ scheme
in \cite{LXKGLG11,Mitzenmacher16,Stolyar15b} all assume that the
loads at the various dispatchers are strictly equal.
We are not aware of any results for heterogeneous dispatcher loads.
To the best of our knowledge, power-of-$d$-policies have not been
considered in a multiple-dispatcher scenario at all.
While the results in~\cite{Stolyar15b} show that the JIQ scheme
is asymptotically optimal for symmetric dispatcher loads, even when
the servers are heterogeneous, it is readily seen that power-of-$d$
policies cannot even be maximally stable in that case for any fixed
value of~$d$. 

{\em Organization of the paper.} 
The remainder of the paper is organized as follows.
In Section~\ref{mode} we present a detailed model description,
specify the two proposed enhancements and state the main results.
In Section~\ref{jacksonnetwork} we describe how the blocking
scenario can be represented in terms of a closed Jackson network,
and leverage the associated product-form distribution to obtain
an insightful formula for the blocking probability.
We then turn to a fluid-limit approach
in Section~\ref{fluidlimitinblocking} to analyze the two proposed
enhancements in the blocking scenario.
A similar analysis is adopted in Section~\ref{fluidlimitinqueueing}
in the queueing scenario to obtain results for the basic model
and the enhanced variants.
Finally, in Section~\ref{conc} we make some concluding remarks
and briefly discuss future research directions.

\section{Model description, notation and key results}\label{mode}
\begin{figure}[t!!]
\begin{center}
\includegraphics[width=\scaling\columnwidth]{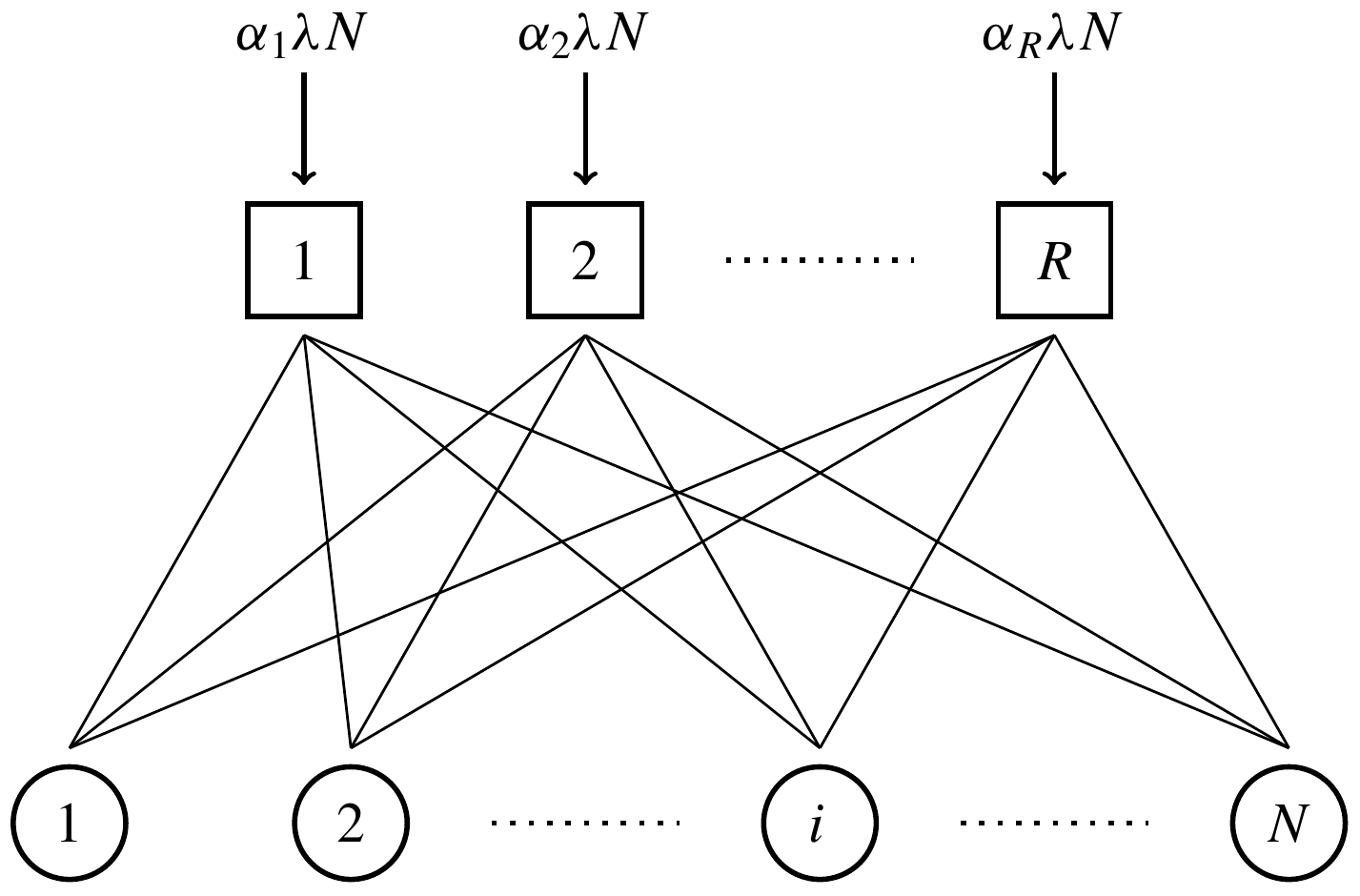}
\caption{Schematic view of the model with $R$~dispatchers and $N$~servers.}
\label{modelfigure}
\end{center}
\end{figure}

We consider a system with $N$~parallel identical servers
and a fixed set of $R$ (not depending on $N$) dispatchers, as depicted in \cref{modelfigure}.
Jobs arrive at dispatcher~$r$ as a Poisson process of rate
$\alpha_r \lambda N$, with $\alpha_r > 0$, $r = 1, \dots, R$,
$\sum_{r = 1}^{R} \alpha_r = 1$, and $\lambda$ denoting the job
arrival rate per server.
For conciseness, we denote $\alpha = (\alpha_1, \dots, \alpha_R)$,
and without loss of generality we assume that the dispatchers are
indexed such that $\alpha_1 \geq \alpha_2 \geq \dots \geq \alpha_R$.
The job processing requirements are independent and exponentially
distributed with unit mean at each of the servers.

When a server becomes idle, it sends a token to one of the
dispatchers selected uniformly at random, advertising its availability.
When a job arrives at a dispatcher which has tokens available,
one of the tokens is selected, and the job is immediately
forwarded to the corresponding server.

We distinguish two scenarios when a job arrives at a dispatcher
which has no tokens available, referred to as the {\em blocking\/}
and {\em queueing\/} scenario respectively.
In the blocking scenario, the incoming job is blocked
and instantly discarded.
In the queueing scenario, the arriving job is forwarded to one
of the servers selected uniformly at random.
If the selected server happens to be idle, then the outstanding
token at one of the other dispatchers is revoked.

In the queueing scenario we assume $\lambda < 1$, which is not only
necessary but also sufficient for stability.
It is not difficult to show that the joint queue length process is
stochastically majorized by a case where each job is sent to
a uniformly at random selected server.
In the latter case, the system decomposes into $N$~independent
M/M/1 queues, each of which has load $\lambda<1$ and is stable.

Denote by $X_0(t)$ the number of busy servers and by $X_r(t)$ the
number of tokens held by dispatcher~$r$ at time~$t$, $r = 1, \dots, R$.
Note that $\sum_{r = 0}^{R} X_r(t) \equiv N$ for all~$t$.
Also, denote by $Y_i(t)$ the number of servers with $i$~jobs
(including a possible job being processed) at time~$t$, $i \geq 0$,
so that $X_0(t) \equiv \sum_{i = 1}^{\infty} Y_i(t)$.

In the blocking scenario, no server can have more than one job,
i.e. $Y_i(t)=0$ for all $i \geq 2$ and $X_0(t) = Y_1(t)$.
Because of the symmetry among the servers, the state of the system can
thus be described by the vector $X(t) = (X_0(t), X_1(t), \dots, X_R(t))$,
and $\{X(t)\}_{t \geq 0}$ evolves as a Markov process, with state space
$S := \{n \in \mathbb{N}^{R+1}: \sum_{i = 0}^{R} n_i = N\}$.

Likewise, in the queueing scenario, the state of the system can
be described by the vector $U(t) = (Y(t), X_1(t), \dots, X_R(t))$
with $Y(t) = (Y_i(t))_{i \geq 0}$, and $\{U(t)\}_{t \geq 0}$
also evolves as a Markov process.
%$$T := \{(m, n) \in \mathbb{N}^{\infty} \times \mathbb{N}^{R}:
%\sum_{i= 0}^{\infty} m_i + \sum_{i = 1}^{R} n_i = N\}.$$

Denote by $B(R, N, \lambda, \alpha)$ the steady-state blocking
probability of an arbitrary job in the blocking scenario.
Also, denote by $W(R, N, \lambda, \alpha)$ a random variable
with the steady-state waiting-time distribution of an arbitrary job
in the queueing scenario.
%REMOVEDLAST
%We include $R$, $N$, $\lambda$ and $\alpha$ in the notation
%in order to reflect the dependence of both performance metrics
%on these system parameters.
%However, we will occasionally suppress (some of) the arguments for
%brevity when these are fixed or otherwise clear from the context. \\

%We consider two different type of models; the blocking model and the waiting model. In both models, we consider $R$ dispatchers, $N$ identical servers and jobs arriving at dispatcher $r$ according to a Poisson process with rate $\alpha_r\lambda N$ with $\sum \alpha_r=1$. Furthermore, we assume that $\alpha_1 \geq \alpha_2 \geq \hdots \geq \alpha_R$. Jobs require unit mean exponential processing time and a server can only host one job at a time. When a server becomes idle, it sends a token to one of the dispatchers. We first assume that this dispatcher is chosen uniformly at random. When a job arrives at a dispatcher which has tokens available, the job uses one of the tokens and is sent to the corresponding server.

%The blocking model assumes that there are no queues in front of the servers. When a job arrives at a dispatcher that lacks tokens, the job gets blocked. We have the following theorem for the blocking probability $B$.

In \cref{jacksonnetwork} we will prove the following theorem
for the blocking scenario.

\begin{mytheorem}[Least-loaded dispatcher determines blocking]
\label{maintheorem}
As $N\to\infty$,
\begin{equation*}
B(R,N,\lambda,\alpha) \to \max\{1-R\alpha_R,1-1/\lambda\}.
\end{equation*}
\end{mytheorem}

\cref{maintheorem} shows that in the many-server limit the system
performance in terms of blocking is either determined by the relative
load of the least-loaded dispatcher, or by the aggregate load.
This may be informally explained as follows.
Let $\Nb$ be the expected fraction of busy servers in steady state,
so that each dispatcher receives tokens on average at a rate $\Nb N/R$.
We distinguish two cases, depending on whether a positive
fraction of the tokens reside at the least-loaded dispatcher~$R$
in the limit or not.
If that is the case, then the job arrival rate $\alpha_R \lambda N$
at dispatcher~$R$ must equal the rate $\Nb N / R$
at which it receives tokens, i.e., $\Nb / R = \alpha_R \lambda$.
Otherwise, the job arrival rate $\alpha_R \lambda N$ at
dispatcher~$R$ must be no less the rate $\Nb N / R$ at which it
receives tokens, i.e., $\Nb / R \leq \alpha_R \lambda$.
Since dispatcher~$R$ is the least-loaded, it then follows that
$\Nb / R \leq \alpha_r \lambda$ for all $r = 1, \dots, R$,
which means that the job arrival rate at all the dispatchers is
higher that the rate at which tokens are received.
Thus the fraction of tokens at each dispatcher is zero in the limit,
i.e., the fraction of idle servers is zero, implying $\Nb = 1$.
Combining the two cases, and observing that $\Nb \leq 1$,
we conclude $\Nb = \min\{R \alpha_R \lambda, 1\}$.
Because of Little's law, $\Nb$ is related to the blocking
probability~$B$ as $\Nb = \lambda (1 - B)$.
This yields $1 - B = \min\{R \alpha_R \lambda, 1 / \lambda\}$,
or equivalently, $B = \max\{1 - R \alpha_R, 1 - 1 / \lambda\}$
as stated in Theorem~1.

The above explanation also reveals that in the limit dispatcher~$R$
(or the set of least-loaded dispatchers in case of ties) inevitably
ends up with all the available tokens, if any.
The accumulation of tokens hampers the visibility of idle servers
to the heavier-loaded dispatchers, and leaves idle servers stranded
while jobs queue up at other servers.
 
\begin{figure}[t!!!]
\begin{center}
\includegraphics[width=\scaling\columnwidth]{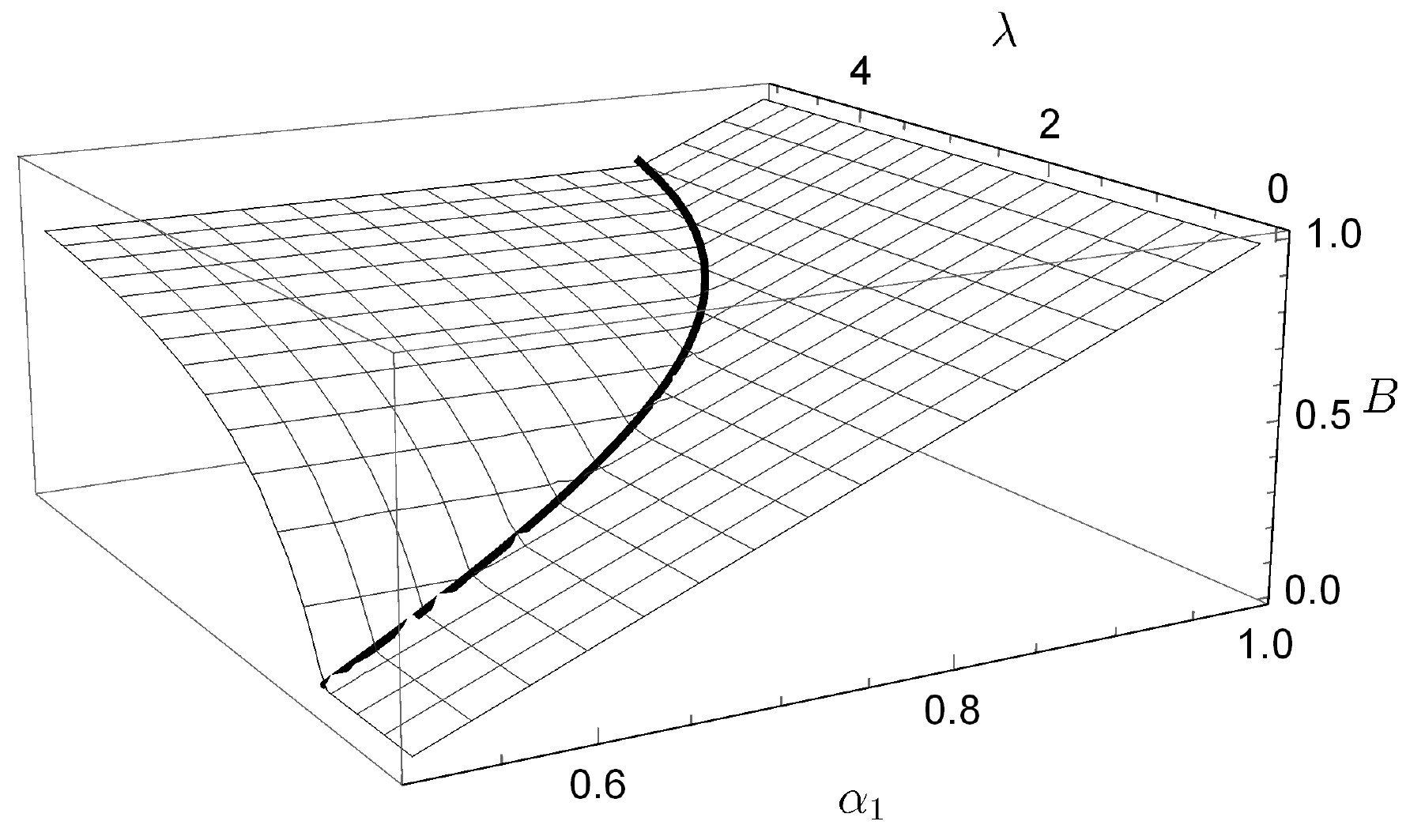}
\caption{Blocking probability $B(2, N, \lambda, (\alpha_1, 1 - \alpha_1))$
obtained by \cref{generalformulaBlocking},
for $R=2$ dispatchers and $N=10^5$ servers as function
of~$\lambda$ and~$\alpha_1$.}
\label{3db1}
\end{center}
\end{figure}

\cref{3db1} illustrates \cref{maintheorem} for $R=2$
dispatchers and $N=10^5$ servers, and clearly reflects the two
separate regions in which the blocking probability depends
on either~$\alpha_R$ or~$\lambda$. The line represents the cross-over curve $R \alpha_R = 2 \alpha_2 = 2 (1 - \alpha_1) = 1 / \lambda$.

In Section~\ref{fluidlimitinqueueing} we will establish the following theorem
for the queueing scenario.

\begin{mytheorem}[Mean waiting time]\label{maintheorem2}
For $\lambda < 1$ and $N\to\infty$,

\[
\mathbb{E}[W(R, N, \lambda, \alpha)] \to
\frac{\lambda_2(R, \lambda, \alpha)}{1 - \lambda_2(R, \lambda, \alpha)},
\]
where 
\[
\lambda_2(R,\lambda,\alpha)=1-\frac{1-\lambda\sum_{i=1}^{r^*}\alpha_i}{1-\lambda r^* / R}
\]
with 
\begin{equation*}\label{maintheorem2r}
r^*=\sup
\big\{r \big| \alpha_r > \frac{1}{R}\frac{1-\lambda\sum_{i=1}^{r}\alpha_i}{1-\lambda r/R}\big\}
\end{equation*}
and the convention that $r^*=0$ if $\alpha_1=\hdots=\alpha_R=1/R$.
\end{mytheorem}
\noindent
$\lambda_2$ can be interpreted as the rate at which jobs are forwarded to randomly selected servers. Furthermore, dispatchers $1, \hdots, r^*$ receive tokens at a lower rate than the incoming jobs, and in particular $\lambda_2^*=0$ if and only if $r^*=0$. When $R=2$, \cref{maintheorem2} simplifies to
\begin{equation*}\label{maintheorem22}
\mathbb{E}[W(2, N, \lambda, (1 - \alpha_2, \alpha_2))] \rightarrow
\frac{\lambda (1-2\alpha_2)}{2-2\lambda (1-\alpha_2)}.
\end{equation*}

When the arrival rates at all dispatchers are strictly equal,
i.e., $\alpha_r = 1/R$ for all $r = 1, \dots, R$, \cref{maintheorem,maintheorem2}
indicate that the stationary blocking probability and the
mean waiting time asymptotically vanish in a regime where the total
number of servers~$N$ grows large, which is in agreement with
the results in~\cite{Stolyar15b}.
However, when the arrival rates at the various dispatchers are not
perfectly equal, so that $\alpha_R < 1 / R$, the blocking probability and mean wait are strictly
positive in the limit, even for arbitrarily low overall load
and an arbitrarily small degree of skewness in the arrival rates.
Thus, the basic JIQ scheme fails to achieve asymptotically optimal
performance when the dispatcher loads are not strictly equal.

In order to counter the above-described performance degradation for
asymmetric dispatcher loads, we propose two enhancements.

\begin{myenhancement}[Non-uniform token allotment]\label{alg1}
When a server becomes idle, it sends a token to dispatcher~$r$
with probability~$\beta_r$.
\end{myenhancement}

\begin{myenhancement}[Token exchange mechanism]\label{alg2}
Any token is transferred to a uniformly randomly selected
dispatcher at rate~$\nu$.
\end{myenhancement}

Note that the token exchange mechanism only creates a constant
communication overhead per job as long as the rate~$\nu$ does not
depend on the number of servers~$N$, and thus preserves the
scalability of the basic JIQ scheme.

The above enhancements can achieve asymptotically optimal
performance for suitable values of the $\beta_r$ parameters
and the exchange rate~$\nu$, as stated in the next proposition.

\begin{myproposition}[Vanishing blocking and waiting]\label{maintheorem3}
For any $\lambda<1$, the stationary blocking probability in the blocking scenario and the mean waiting time in the queueing scenario asymptotically vanish as $N\to\infty$, upon using \cref{alg1} with $\beta_r=\alpha_r$ or \cref{alg2} with $\nu\geq \frac{\lambda}{1-\lambda}(\alpha_1 R-1)$.
\end{myproposition}

%\begin{myproposition}[Vanishing wait]\label{maintheorem4}
%The mean waiting time in the queueing scenario asymptotically vanishes as $N\to\infty$, upon using \cref{alg1} with $\beta_r=\alpha_r$ or \cref{alg2} with $\nu\geq \frac{\lambda}{1-\lambda}(\alpha_1 R-1)$.
%\end{myproposition}

The minimum value of $\nu$ required in the blocking scenario may be intuitively understood as follows. Zero blocking means that a fraction $\lambda$ of the servers must be busy, and thus a fraction $1-\lambda$ of the tokens reside with the various dispatchers, while the heaviest loaded dispatcher 1 receives enough tokens for all incoming jobs: $\alpha_1  \lambda \leq \lambda/R + \nu(1-\lambda)/R$ which is satisfied by the given minimum value of $\nu$.

%\colorbox{green}{We} will intuitively explain the minimal value of $\nu$ in the blocking scenario. \colorbox{green}{Observe} the bottleneck of the system; the heaviest-loaded dispatcher 1. Jobs arrive with rate $\alpha_1\lambda$ and it will have (close to) 0 tokens. Suppose that a fraction $x$ of the servers is busy. Then, tokens arrive with rate $x/R$. Since $x$ servers are busy, $1-\lambda$ tokens are at the dispatchers. With \cref{alg2} in place, another $(1-\lambda)\nu/R$ tokens will arrive at dispatcher 1. Since no jobs can be blocked, the rate of incoming tokens should be at least the rate of arriving jobs. This leads to $\alpha_1\lambda \leq x/R + \nu \frac{1-x}{R}$, which is satisfied by the given minimal value of $\nu$.

A similar reasoning applies to the queueing scenario, although in that case the number of servers with exactly one job no longer equals the number of busy servers, and a different approach is needed.

In order to establish \cref{maintheorem2,maintheorem3}, we examine 
in \cref{fluidlimitinblocking,fluidlimitinqueueing} the fluid limits for the blocking and queueing
scenarios, respectively. 
Rigorous proofs to establish weak convergence to the fluid limit
are omitted, but can be constructed
along similar lines as in~\cite{HK94}.
The fluid-limit regime not only provides mathematical tractability,
but is also particularly relevant given the massive numbers
of servers in data centers and cloud operations.
Simulation experiments will be conducted to verify the accuracy
of the fluid-limit approximations, and show an excellent match,
even in small systems (small values of $N$).

\section{Jackson network representation}
\label{jacksonnetwork}

%We consider $R$ dispatchers, $N$ identical servers and jobs arriving to dispatcher $r$ according to a Poisson process with rate $\alpha_r \lambda N$ with $\sum \alpha_r=1$. Furthermore, we assume that $\alpha_1 \geq \alpha_2 \geq \hdots \geq \alpha_R$.
%The blocking model assumes that servers can only host one job, which means that no queues can be formed in front of the servers. Jobs require unit mean exponential processing times. When a server becomes idle, it sends a token to one of the dispatchers. We assume that this dispatcher is chosen uniformly at random.
%Under these assumptions, jobs can be blocked from the system. Upon arrival to a dispatcher, a job is blocked when the dispatcher has no tokens available. In case there are tokens, the job is transferred to one of the idle servers and the corresponding token is removed from the dispatcher. Here we assume that tokens are used in a FCFS manner.
In this section we describe how the blocking scenario can be
represented in terms of a closed Jackson network.
We leverage the associated product-form distribution to express
the asymptotic blocking probability as a function of the aggregate
load and the minimum load across all dispatchers,
proving \cref{maintheorem}.

We view the system dynamics in the blocking scenario in terms
of the process $\{X(t)\}_{t \geq 0}$ as a fixed total population
of $N$~tokens that circulate through a network of $R+1$ stations.
Specifically, the tokens can reside either at station~$0$,
meaning that the corresponding server is busy, or at some station~$r$,
indicating that the corresponding server is idle and has
an outstanding token with dispatcher~$r$, $r = 1, \dots, R$.
%The system consists of the $R$~dispatchers and the $N$~servers,
%but we only need to keep track of the number $X_0(t)$ of busy
%servers (that have a token) and the number of tokens $X_r(t)$
%at each dispatcher~$r$, $r = 1, \dots, R$.
%that circulate through the system.

Let $s_i(k)$ denote the service rate at station~$i$ when there are
$k$~tokens present.
Then $s_0(k)=k$ and $s_r(k)=1$ for $r=1,\dots,R$.
The service times are exponentially distributed at all stations,
but station~$0$ is an infinite-server node with mean service
time $\mu_0^{-1}=1$, while station~$r$ is a single-server node with
mean service time $\mu_r^{-1} = (\alpha_r \lambda N)^{-1}$, $r=1,\dots,R$.
The routing probabilities $p_{ij}$ of tokens moving from station~$i$
to station~$j$ are given by $p_{r0}=1$ for $r=1,\dots,R$
and $p_{0r}=1/R$ for $r=1,\dots,R$.
With $\gamma_i$ denoting the throughput of tokens at station~$i$,
the traffic equations
\[
\gamma_i = \sum_{j = 0}^{R} \gamma_j p_{ji}, \quad i=0,\dots,R
\]
uniquely determine the relative values of the throughputs.
%REMOVEDLAST
%i.e., the throughputs up to a common scaling factor.

%The probabilities~$\pi(n)$ satisfy the balance equations
%\[
%\left[ \sum\limits_{i = 0}^{R} \mu_i r_i(n_i) \right] \pi(n) =
%\sum\limits_{i = 0}^{R} \sum\limits_{j = 1}^{M}
%\mu_i r_i(n_i + 1) p_{ij} \id{n_j > 0} \pi(n + e_i - e_j)
%\]
%for all states $n \in S$, along with the normalization condition
%\[
%\sum\limits_{n \in S} \pi(n) = 1.
%\]

Let $\pi(n) := \lim_{t \to \infty} \mathbb{P}(X(t) = n)$ be the
stationary probability that the process $\{X(t)\}_{t \geq 0}$
resides in state $n \in S$.
The theory of closed Jackson networks~\cite{Kelly79} implies
\[
\pi(n_0,n_1,\hdots,n_R)=G^{-1} \prod_{i=0}^R  \frac{(\gamma_i/\mu_i)^{n_i}}{\prod_{m=1}^{n_i} s_i(m)},
\]
with $G$ a normalization constant.

The blocking probability can then be expressed by summing the
probabilities $\pi(n)$ over all the states with $n_r = 0$ where no
tokens are available at dispatcher~$r$, and weighting these
with the fractions $\alpha_r$, $r = 1, \dots, R$:
%We then consider the probability that a job arrives at some dispatcher with no tokens available, so that the job gets blocked from the system. We denote this blocking probability by $B(R,N,\lambda,\alpha)$, and note that it can be expressed in terms of the stationary distribution $\pi(n)$ by summing over all states... KIJKEN HOE DIT GOED TE BESCHRIJVEN.
\begin{equation}\label{generalformulaBlocking}
B(R,N,\lambda,\alpha)=\sum_{r=1}^R \alpha_r \sum_{n\in \{n | n_r=0\}} \pi(n).
\end{equation}
Despite this rather complicated expression,
\cref{maintheorem} provides a compact characterization of the
blocking probability in the many-server limit $N\to \infty$,
as will be proved in the Appendix.
The proof uses stochastic coupling, for which we define a `better' system and a `worse' system. Both systems are amenable to analysis and have an identical blocking probability in the many-server limit $N\to \infty$.

The better system merges the first $R-1$ dispatchers into one super-dispatcher, which results in two dispatchers with arrival rates $(1-\alpha_R)\lambda N$ and $\alpha_R\lambda N$, respectively. However, in contrast to the original blocking scenario, when a job is completed and leaves a server idle, a token is not sent to either dispatcher with equal probability. Instead, tokens are sent to the super-dispatcher with probability $\frac{R-1}{R}$. To analyze this better system, we study in \cref{modelA} the blocking scenario with $R = 2$ enhanced with non-uniform token allotment.
%and the additional feature of a heterogeneous distribution of tokens over the dispatchers. This model is called the ``Heterogeneous blocking model with two dispatchers'' and analyzed in \cref{modelA}.

The worse system thins the incoming rates of jobs at the dispatchers,
so that some jobs are blocked, irrespective of whether or not the
dispatcher has any tokens available.
This thinning process is defined as follows:
a job arriving at dispatcher~$r$ is blocked with probability
\[%ADDEDMARK287
\frac{\alpha_r-\alpha_R}{\alpha_r} +
\frac{\max\left\{0, \alpha_R-1/(R\lambda)\right\}}{\alpha_r}.
\]
This thinning process is designed in such a way that the system with admitted jobs behaves as a system with total arrival rate $\lambda\leq 1$ in which all arrival rates are equal ($\alpha_1=\hdots=\alpha_R$), which is analyzed in \cref{modelB}.
% To analyze the worse system, we study the ``Symmetric blocking model with arbitrarily many dispatchers'' in \cref{modelB}.

%ADDEDMARK287
With coupling, one can show that the blocking probability of the
`better system' is lower and of the `worse system' is higher,
which completes the proof.
Specifically, when the arrival moments, the service times and the
token-allotment are coupled, the number of tokens used at each
dispatcher by time~$t$ is always lower in the worse system
and higher in the better system.
Due to page limitations, the detailed coupling arguments are omitted,
but it is intuitively clear that the better system performs better
and the worse system performs worse.
Namely, the tokens at dispatchers 1 to $R-1$ are consolidated in the better system. If there is at least one token amongst these dispatchers, any job arriving at any of the dispatchers can make use of a token. In the original system, a job is blocked when the token amongst the first $R-1$ dispatchers, is not present at the dispatchers at which a job arrives. The worse system performs obviously worse, since blocking jobs beforehand has no benefits for the acceptance of jobs.

While we assumed exponentially distributed service times, the infinite-server node is symmetric and thus the product-form solution in \cref{generalformulaBlocking} as well as \cref{maintheorem} still hold for phase-type distributions \cite{Kelly79}.

\section{Fluid limit in the blocking scenario}
\label{fluidlimitinblocking}

We now turn to the fluid-limit analysis and start with the blocking
scenario.
We consider a sequence of systems indexed by the total number of
servers~$N$.
Denote by $x_0^N(t) = \frac{1}{N} X_0^N(t)$ the fraction of busy
servers and by $x_r^N(t) = \frac{1}{N} X_r^N(t)$ the normalized number of tokens held by dispatcher $r$ in the $N$-th system at time~$t$.
Further define $x^N(t) = (x_0^N(t), x_1^N(t) \dots, x_R^N(t))$
and assume that $x^N(0) \to x^\infty$ as $N \to \infty$,
with $\sum_{i = 0}^{R} x_i^\infty = 1$.
Then any weak limit $x(t)$ of the sequence $\{x^N(t)\}_{t \geq 0}$
as $N \to \infty$ is called a fluid limit. 

The fluid limit $x(t)$ in the blocking scenario
with \cref{alg1,alg2} in place satisfies the set of differential
equations
\begin{equation}\label{fpexpr1}
\frac{\dd x_0(t)}{\dd t} = \sumR z_r(t) - x_0(t),
\end{equation}
\begin{equation}\label{fpexpr2}
\frac{\dd x_r(t)}{\dd t} =
\beta_r x_0(t) + \nu\left(\frac{1-x_0(t)}{R}-x_r(t)\right) - z_r(t),
\end{equation}
with
\begin{equation}\label{fpexpr2z}
z_r(t) = \alpha_r \lambda -
\left[\alpha_r \lambda - \beta_r x_0(t) - \nu \frac{1-x_0(t)}{R}\right]^+
\mathds{1}\left\{x_r(t)=0\right\},
\end{equation}
where $[\cdot]^+=\max\{\cdot,0\}$ and initial condition $x(0) = x^\infty$.

%{\tt Next sentence optional:}
%Note that $\sum_{i = 0}^{R} \frac{\dd x_i(t)}{\dd t} =
%1 - \sum_{i = 0}^{R} x_i(t)$, so $\sum_{i = 0}^{R} x_i^\infty = 1$
%implies $\sum_{i = 0}^{R} x_i(t) \equiv 1$ for all $t \geq 0$.

The above set of fluid-limit equations may be interpreted as follows.
The term $z_r(t)$ represents the (scaled) rate at which dispatcher~$r$
uses tokens and forwards incoming jobs to idle servers at time~$t$.
\cref{fpexpr2z} reflects that the latter rate equals the job
arrival rate $\alpha_r \lambda$, unless the fraction of tokens held
by dispatcher~$r$ is zero ($x_r(t) = 0$), and the rate
$\beta_r x_0(t) + \nu (1 - x_0(t)) / R$ at which it receives
tokens from idle servers or through the exchange mechanism is less
than the job arrival rate. \cref{fpexpr1} states that the rate of change in the fraction
of busy servers is the difference between the aggregate rate
$\sumR z_r(t)$ at which the various dispatchers use tokens
and forward jobs to idle servers, and the rate $x_0(t)$
at which jobs are completed and busy servers become idle.
%The total usage rate of tokens; the rate at which servers become busy, is represented in the first term of \cref{fpexpr1}. The second term is the scaled rate at which servers become idle, which is $x_0(t)$ because all jobs require a unit mean time to be served.
\cref{fpexpr2} captures that the rate of change of the fraction
of tokens held by dispatcher~$r$ is the balance of the rate
$\beta_r x_0(t) + \nu (1 - x_0(t)) / R$ at which it receives tokens
from idle servers or through the exchange mechanism,
and the rate $z_r(t) + \nu x_r(t)$ at which it uses tokens
and forwards jobs to idle servers or releases tokens through the
exchange mechanism. \\

%We turn to \cref{fpexpr2}, in which the first term is the scaled rate at which servers become idle, multiplied by the probability that the token is sent to dispatcher $r$. A token is used at dispatcher $r$ with rate $z_r(t)$. The final term represents \cref{alg2}. With rate $\nu$, a token is moved. The number of outstanding tokens is $(1-x_0(t))$, and with probability $1/R$, a token is sent to dispatcher $r$ with rate $\nu$. However, one of the $x_r(t)$ at dispatcher $r$ also leave with rate $\nu$. This explains the final term of \cref{fpexpr2}.

\begin{figure}[t!!!]
	\begin{center}
		\includegraphics[width=\scaling\columnwidth]{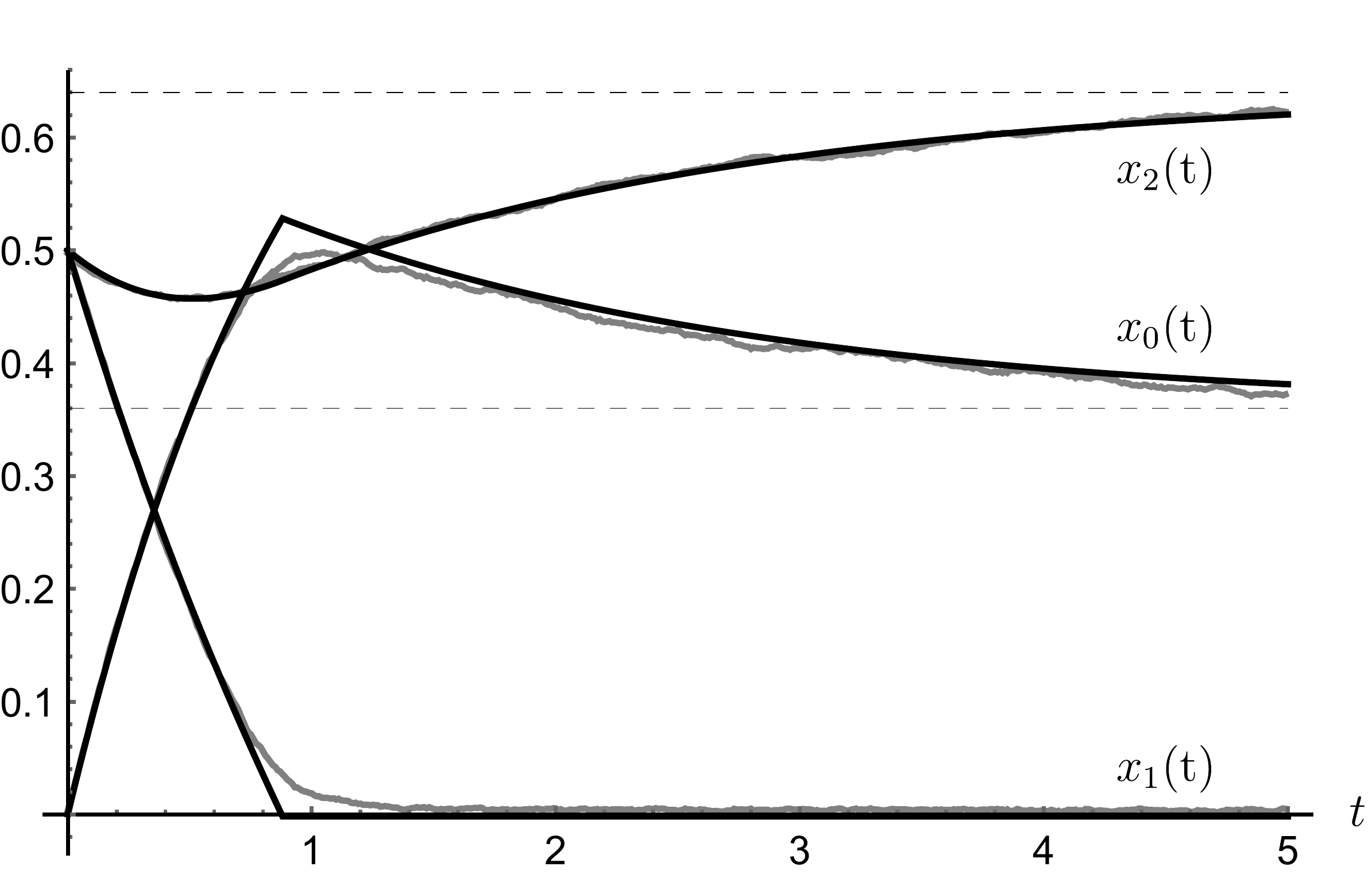}
		\caption{Fluid-limit trajectories $x_i(t)$ for $R=2$ dispatchers,
			$\lambda=0.9$ and $\alpha_1=0.8$.
			Averaged sample paths for $N=100$ servers are represented by lighter
			colors, and closely match the fluid-limit dynamics.}\label{fwtp5}
	\end{center}
\end{figure}

\cref{fwtp5} shows the exact and simulated fluid-limit trajectories.
We observe that the simulation results closely match the
fluid-limit dynamics.
We further note that in the long run only dispatcher~$2$ with the
lower arrival rate holds a strictly positive fraction of the tokens,
corroborating Theorem~1.

\subsection{Fixed-point analysis}

In order to determine the fixed point(s)~$x^*$, we set
$\dd x_i(t)/\dd t = 0$ for all $i = 0, 1, \dots, R$, and obtain
\begin{equation}\label{fpexpr3c}
x_0^* = \sumR z_r^*,
\end{equation}
\begin{equation}\label{fpexpr3d}
z_r^* = \beta_r x_0^* + \nu\left(\frac{1-x_0^*}{R}-x_r^*\right),
\end{equation}
and
\begin{equation}
\label{fpexpr3b}
z_r^* = \alpha_r \lambda -
\left[\alpha_r \lambda - \beta_r x_0^* - \nu \frac{1-x_0^*}{R}\right]^+
\mathds{1}\left\{x_r^*=0\right\}.
\end{equation}

Without proof, we assume that the many-server ($N\to\infty$) and stationary ($t\to\infty$) limits commute, so that $x_0^*$ is also the limit of the mean fraction of busy servers in stationarity. Because of Little's law, the limit $B$ of the blocking probability satisfies 
\begin{equation}\label{littlecons}
x_0^*=\lambda(1-B).
\end{equation}
This in particular implies that $x_0^*=\lambda$ leads to $B=0$: vanishing blocking. 

\emph{Basic JIQ scheme.}
We first consider the basic JIQ scheme, i.e., $\beta_r = 1 / R$
for all $r = 1, \dots, R$ and $\nu = 0$.
\cref{fpexpr3d,fpexpr3b} yield
\begin{equation*}
\frac{x_0^*}{R} = z_r^* = \alpha_r \lambda -
\left[\alpha_r \lambda - \frac{x_0^*}{R}\right]^+
\mathds{1}\left\{x_r^*=0\right\},
\end{equation*}
or equivalently,
\begin{equation}
\label{basicjiq}
\alpha_r \lambda - \frac{x_0^*}{R} = 
\left[\alpha_r \lambda - \frac{x_0^*}{R}\right]^+
\mathds{1}\left\{x_r^*=0\right\}.
\end{equation}
Now let ${\mathcal I} = \{r: \alpha_r = \alpha_R\}$ be the index
set of the least-loaded dispatchers.
\cref{basicjiq} forces $x_r^* = 0$ for all
$r \notin {\mathcal I}$.

We now distinguish two cases, depending on whether or not
$x_r^* = 0$ for all $r \in {\mathcal I}$ as well.
If that is the case, then we must have $x_0^* = 1$,
and $\alpha_R \lambda \geq x_0^*/R$, i.e.,
$\lambda \geq 1 / (R \alpha_R)$.
Otherwise, we must have $\alpha_R \lambda = x_0^*/R$,
i.e., $x_0^* = R \alpha_R \lambda$, so $x_0^* \leq 1$ forces
$\lambda \leq 1 / (R \alpha_R)$.

In conclusion, we have $x_0^* = \min\{R \alpha_R \lambda, 1\}$.
When $\lambda \geq 1 / (R \alpha_R)$ so that $x_0^* = 1$,
it must be the case that $x_r^* = 0$ for all $r = 1, \dots, R$.
When $\lambda < 1 / (R \alpha_R)$ so that $x_0^* < 1$,
any vector~$x^*$ with $x_r^* = 0$ for all $r \notin {\mathcal I}$ and $\sum_{r \in {\mathcal I}} x_r^* = 1 - x_0^*$ is a fixed point.
In particular, for equal dispatcher loads, i.e., $\alpha_r = 1 / R$
for all $r = 1, \dots, R$, so that ${\mathcal I} = \{1, \dots, R\}$,
we have $x_r^* = 0$ for all $r = 1, \dots, R$ when $\lambda \geq 1$,
while any vector~$x^*$ with $\sumR x_r^* = 1 - \lambda$ is a fixed
point when $\lambda < 1$. 

We use \cref{littlecons} to find
$B = 1 - \frac{1}{\lambda} \min\{R \alpha_R \lambda, 1\} =
\max\{1 - R \alpha_R, 1 - 1 / \lambda\}$, which agrees with \cref{maintheorem}.

In \cref{table1} we compare the fluid-limit approximations
for the blocking probability with the exact formula from the
Jackson network representation and simulation results for various
numbers of servers.
\begin{table}[h]
\begin{center}
\caption{Blocking probabilities for $\lambda=0.9$ and $\alpha_1=0.8$ or $\alpha_1=0.6$.}
\begin{tabular}{ r  c  c  c  c }
\hline
& \multicolumn{2}{c}{$\lambda=0.9, \alpha_1=0.8$} & \multicolumn{2}{c}{$\lambda=0.9, \alpha_1=0.6$}\\
N & Jackson & simulation & Jackson & simulation\\
\midrule
$10$    & 0.6021 & 0.6032 & 0.3201 & 0.3205 \\
$20$    & 0.6000 & 0.6006 & 0.2545 & 0.2552 \\
$50$    & 0.6000 & 0.6004 & 0.2092 & 0.2095 \\
$100$   & 0.6000 & 0.6007 & 0.2006 & 0.2010  \\
%$200$    & 0.6000 & 0. & 0.2000 &  \\
%$500$    & 0.6000 & 0. & 0.2000 &   \\
%$1000$ & 0.6000 & 0.6008 & 0.2000 &   \\
\hline
fluid ($N=\infty$) & 0.6000 & - & 0.2000 & -\\
\hline
\end{tabular}
\label{table1}
\end{center}
\end{table}
%MARKADDED297
\cref{table1} shows that the Jackson network analysis agrees
with the simulation results.
Furthermore, the more symmetric the loads, the lower the blocking
probability, which is consistent with \cref{maintheorem}.
Also, the fluid-limit approximation is highly accurate,
even for a fairly small number of servers.

\subsection{Enhancements}

We now examine the behavior of the system for \cref{alg1,alg2},
and show that they can achieve asymptotically zero blocking
for any $\lambda\leq 1$ and suitable parameter values as identified in \cref{maintheorem3}. In light of \cref{littlecons} it suffices to show that $x_0^*=\lambda$ for both enhancements.

Consider \cref{alg1}; $\beta_r=\alpha_r$ and $\nu=0$. 
\cref{fpexpr3b,fpexpr3d} give $\lambda-x_0^*=[\lambda-x_0^*]^+\mathds{1}\{x_r^*=0\}$ for all $r$ (this shows $x_0^*\leq \lambda$). Assume that $x_0^*<\lambda$. Then, $x_r^*=0$ for all $r$, which results in $x_0^*=1$, since $\sum_{i=0}^R x_i^*=1$. This contradicts $x_0^*\leq \lambda<1$, so that $x_0^*=\lambda$.

Notice that any point for which $x_0^*=\lambda$ and for which $\sum_{i=0}^R x_i^*=1$, results in $z_r^*=\alpha_r\lambda$ which is in accordance with \cref{fpexpr3c,fpexpr3d,fpexpr3b}, and therefore is a fixed point.\\

We next consider \cref{alg2} with $\beta_r=1/R$ and $\nu\geq \frac{\lambda}{1-\lambda}(R\alpha_1-1)$.
We use \cref{fpexpr3c,fpexpr3b} twice, first they give $x_0^*\leq \lambda$, so that $x_0^*=\lambda-\ve$ for some $\ve\geq0$.
Second, $z_r^*=\alpha_r\lambda$ if the term in brackets in \cref{fpexpr3b} is non-positive. Otherwise,
\[
\begin{split}
z_r^* &\geq \alpha_r\lambda - \left[\alpha_r \lambda - \frac{x_0^*}{R}-\nu\frac{1-x_0^*}{R}\right]\\
&\geq \alpha_r\lambda - \alpha_r\lambda + \frac{\lambda-\ve}{R} + \frac{\lambda}{1-\lambda}(R\alpha_1-1)
\frac{1-\lambda}{R}+\nu\frac{\ve}{R}\\
&\geq \alpha_r \lambda - \frac{\ve}{R}(1-\nu),
\end{split}
\]
which by \cref{fpexpr3c} gives $\lambda-\ve=x_0^*=\sumR z_r^*\geq \lambda-\ve(1-\nu)$, so that $x_0^*=\lambda$.\\

\begin{figure}[t!!!]
\begin{center}
\includegraphics[width=\scaling\columnwidth]{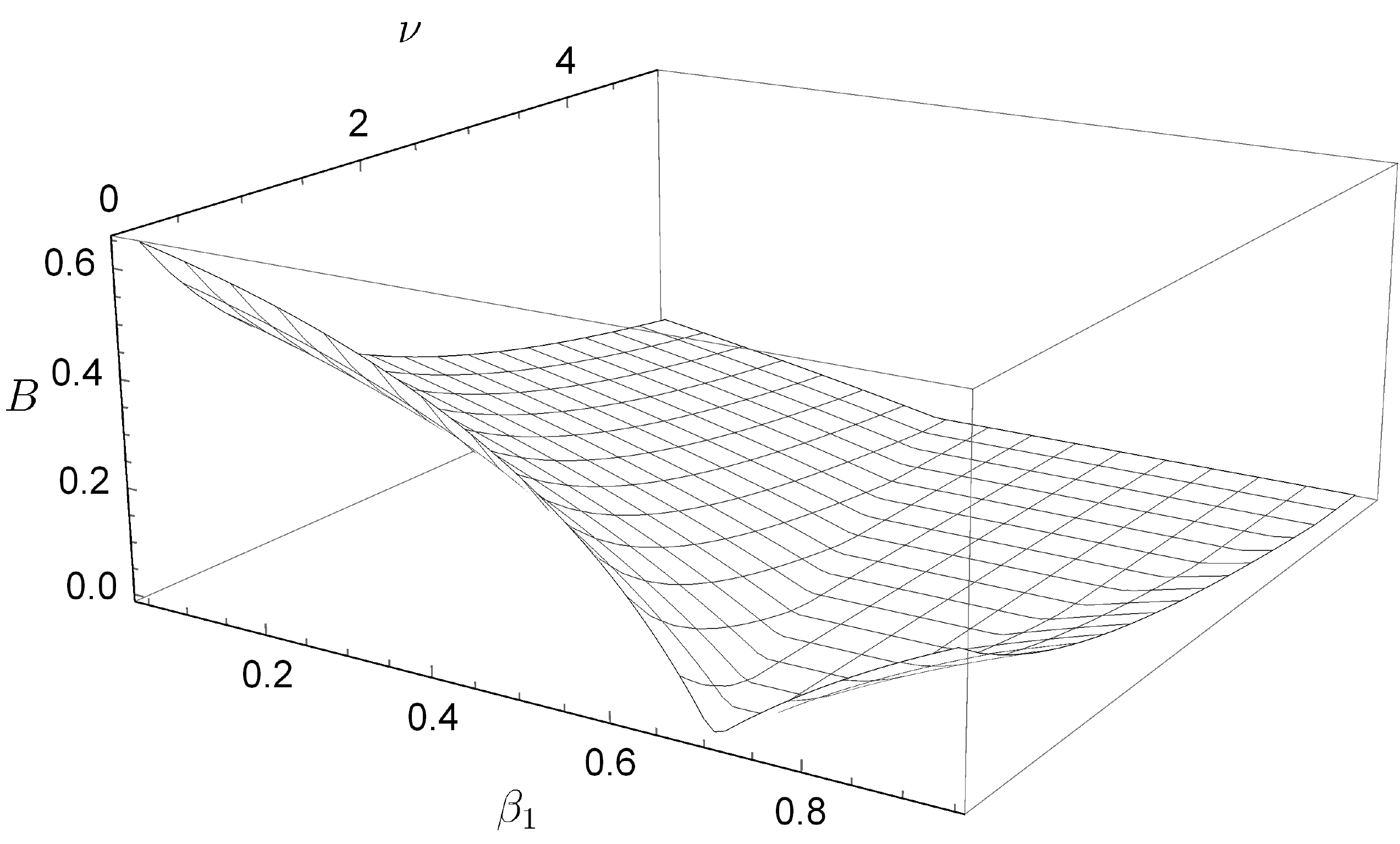}
\caption{Blocking probability $B$ in the limit for $R=2$, $\lambda=0.9$ and $\alpha_1=0.7$, for different values of $\beta_1$ and $\nu$.}\label{fwtp3}
\end{center}
\end{figure}

\cref{fwtp3} displays the blocking probability as $N\to\infty$ for the system with both enhancements.
%ADDEDMARK287
Since $\alpha_1=0.7$, we have that $\beta_1=0.7$ is optimal.
The blocking probability decreases as $\beta_1$ approaches $\alpha_1$
and as $\nu$ increases.
For $\nu>0$, it suffices to choose $\beta_1$ close to $\alpha_1$,
which implies that it is not necessary to know the exact loads,
for the enhancements to be effective.

%As we can see in this expression, the blocking probability only depends on dispatcher $R$; the least-loaded dispatcher. To verify this result, we will analyze the system in equilibrium, for a general number of servers $N$.

\section{Fluid limit in the queueing scenario}
\label{fluidlimitinqueueing}

We now proceed to the queueing scenario (with $\lambda < 1$ for
stability).
As before, we consider a sequence of systems indexed by the total
number of servers~$N$.
Denote by $y_i^N(t) = \frac{1}{N} Y_i^N(t)$ the fraction of servers
with $i$~jobs and by $x_r^N(t) = \frac{1}{N} X_r^N(t)$ the normalized number
of tokens held by dispatcher $r$ in the $N$-th system at time~$t$,
$r = 1, \dots, R$.
Further define $u^N(t) = (y^N(t), x_1^N(t) \dots, x_R^N(t))$,
with $y^N(t) = (y_i^N(t))_{i \geq 0}$, and assume that
$u^N(0) \to u^\infty$ as $N \to \infty$, with
$\sum_{i = 1}^{\infty} y_i^\infty + \sum_{r = 1}^{R} x_r^\infty = 1$.
Then any weak limit $u(t)$ of the sequence $\{u^N(t)\}_{t \geq 0}$
as $N \to \infty$ is called a fluid limit.

The fluid limit $u(t)$ in the queueing scenario
with \cref{alg1,alg2} in place obeys the set of differential equations
\begin{equation}\label{fluidlimitexpr1}
\frac{\dd y_0(t)}{\dd t} =
y_1(t) - \lambda_1(t) - \lambda_2(t) y_0(t),
\end{equation}
\begin{equation}\label{fluidlimitexpr2}
\begin{split}
\frac{\dd y_i(t)}{\dd t} &= \lambda_1(t)\mathds{1}\left\{i=1\right\} + \lambda_2(t) y_{i-1}(t) \\
&+ y_{i+1}(t) - y_i(t) - \lambda_2(t) y_i(t)
\mbox{ for all } i \geq 1,
\end{split}
\end{equation}
\begin{equation}\label{fluidlimitexpr3}
\frac{\dd x_r(t)}{\dd t} =
\beta_r y_1(t) + \nu\left(\frac{y_0(t)}{R}-x_r(t)\right) - z_r(t) -
\lambda_2(t) x_r(t),
\end{equation}
with
\[%^M
z_r(t) = \alpha_r \lambda -
\left[\alpha_r \lambda - \beta_r y_1(t) - \nu \frac{y_0(t)}{R}\right]^+
\mathds{1}\left\{x_r(t)=0\right\},
\]%^M
\begin{equation}\label{fluidlimitexpr4}
\lambda_1(t) = \sumR z_r(t), \quad
\lambda_2(t) = \lambda - \lambda_1(t),
\end{equation}
and initial condition $u(0) = u^\infty$.

%{\tt Next sentence is optional:}
%Note that $\sum_{i = 0}^{\infty} \frac{\dd y_i(t)}{\dd t} +
%\sum_{r = 1}^{R} \frac{\dd x_r(t)}{\dd t} = \sum_{i = 0}^{R} x_i(t) =
%(\lambda_2(t) + \nu) \left(y_0(t) - \sum_{r = 1}^{R} x_r(t)\right) =
%(\lambda_2(t) + \nu) \left(1 - \sum_{i = 0}^{R} x_i(t)\right)$,
%so $\sum_{i = 0}^{\infty} y_i^\infty + \sum_{r = 1}^{R} x_r^\infty = 1$,
%implies
%$\sum_{i = 0}^{\infty} y_i(t) + \sum_{r = 1}^{R} x_r(t) = 1$,
%for all $t \geq 0$.

The above set of fluid-limit equations may be interpreted as follows.
Similarly as in the blocking scenario, the term $z_r(t)$ represents
the (scaled) rate at which dispatcher~$r$ uses tokens and forwards
incoming jobs to idle servers at time~$t$.
Accordingly, $\lambda_1(t)$ is the aggregate rate at which dispatchers
use tokens to forward jobs to (guaranteed) idle servers at time~$t$,
while $\lambda_2(t)$ is the aggregate rate at which jobs are forwarded
to randomly selected servers (which may or may not be idle).
%We define the total usage rate over all dispatchers as $\lambda_1(t)$, which leaves $\lambda_2(t)$ to be the rate at which jobs are dispatched without using a token.
\cref{fluidlimitexpr1} reflects that the rate of change in the
fraction of idle servers is the difference between the aggregate
rate $y_0(t)$ at which jobs are completed by servers with one job,
and the rate $\lambda_1(t)$ at which dispatchers use tokens to forward
jobs to idle servers plus the rate $\lambda_2(t) y_0(t)$ at which
jobs are forward to randomly selected servers that happen to be idle.
%First, we cover \cref{fluidlimitexpr1}. Idle servers can switch state if a type-1 job is allocated to the server (with rate $\lambda_1 N$), when a job is sent to an idle server without any pull message (rate $\lambda_2 y_0 N$). Furthermore, a server can become idle with rate $y_1$.
\cref{fluidlimitexpr2} states that the rate of change in the
fraction of servers with $i$~jobs is the balance of the rate
$\lambda_2(t) y_{i-1}(t)$ at which jobs are forwarded to randomly
selected servers with $i - 1$ jobs plus the aggregate rate
$y_{i+1}(t)$ at which jobs are completed by servers with $i + 1$ jobs,
and the rate $\lambda_2(t) y_i(t)$ at which jobs are forwarded
to randomly selected servers with $i$~jobs plus the aggregate rate
$y_i(t)$ at which jobs are completed by servers with~$i$ jobs.
In case $i = 1$, the rate at which dispatchers use tokens to forward
jobs to idle servers should be included as additional positive term.

%We move to \cref{fluidlimitexpr2} for $i=1$; servers with 1 job switch state with rate $y_1 N$ (job served) and rate $\lambda_2 y_1 N$ (job is sent to the server without pull message). A server can arrive at the 1-job-state with rate $y_2 N$ (a server with 2 jobs finishes service of one of them), rate $\lambda_2(1-x) N$ when a job is sent to an idle server without pull message and rate $\lambda_1 N$ when a job uses a pull message. We obtain \cref{fluidlimitexpr2} for general $i$ by removing the first term, since tokens only dispatch jobs at idle servers. 

\cref{fluidlimitexpr3} is similar to \cref{fpexpr2}, where the
additional fourth term captures the rate at which tokens are revoked
when jobs are forwarded to randomly selected servers that happen
to be idle.

\subsection{Fixed-point analysis}

In order to determine the fixed point(s)~$u^*$, we set
$\dd y_i(t)/\dd t = 0$ for all $i \geq 0$,
and $\dd x_r(t)/\dd t = 0$ for all $r = 1, \dots, R$.
We obtain
\begin{equation}
\label{mm}
y_1^* = \lambda_1^* + \lambda_2^* y_0^*, \hspace{.2cm} (1+\lambda_2^*) y_i^* = \lambda_2^* y_{i-1}^* + y_{i+1}^*
\text{ for all } i \geq 2,
\end{equation}
Solving \cref{mm} gives
\begin{equation*}
y_0^*=1-\lambda, \hspace{.2cm} y_k^* = \lambda (1 - \lambda_2^*) \left(\lambda_2^*\right)^{k-1}
\text{ for all } k \geq 1.
\end{equation*}
Thus the mean number of jobs at a server is
\[
\sum_{k = 1}^{\infty} k y_k^* = \sum_{k=1}^{\infty} k \lambda (1 - \lambda_2^*)
\left(\lambda_2^*\right)^{k-1} = \frac{\lambda}{1-\lambda_2^*}.
\]

As for the blocking scenario, we assume that the many-server and stationary limits commute. Little's law then gives
\begin{equation}\label{littlequeueing}
\sum_{k=1}^\infty k y_k^*=\lambda(\mathbb{E}[W]+1),
\end{equation}
where the left-hand side represents the mean number of jobs at a server in the many-server limit. We use \cref{littlequeueing} to obtain
\begin{equation}\label{meanwaitingtime}
\mathbb{E}[W] = \frac{\sum_{k = 1}^{\infty} k y_k^*}{\lambda} - 1 = \frac{\lambda_2^*}{1 - \lambda_2^*},
\end{equation}
which shows vanishing wait in case $\lambda_2^*=0$.
%{\tt Previous or next argument is optional:}
%The limit of the mean waiting time may also be derived by noting
%that a fraction $\lambda_1^* / \lambda$ of the jobs is forwarded
%to an idle server and incurs zero wait, while a fraction
%$\lambda_2^* / \lambda$ of the jobs is assigned to a randomly
%selected server, so that $\mathbb{E}[W] =
%\frac{\lambda_2^*}{\lambda} m^* = \lambda_2^* / (1 - \lambda_2^*)$,
%confirming the above formula.

We also obtain the following equations for the fixed  point:
\begin{equation}\label{zrexpr}
z_r^* = \beta_r \lambda (1 - \lambda_2^*) - \lambda_2^* x_r^* +
\nu\left(\frac{1-\lambda}{R}-x_r^*\right),
\end{equation}
\begin{equation}\label{zrexpr2}
z_r^* = \alpha_r \lambda - \left[\alpha_r \lambda -
\beta_r \lambda (1 - \lambda_2^*) - \nu \frac{1 - \lambda}{R}\right]^+
\mathds{1}\left\{x_r^*=0\right\},
\end{equation}
and
\begin{equation}\label{zrexprl1}
\lambda_1^* = \sum_{r=1}^R z_r^*.
\end{equation}
We define 
\begin{equation}\label{fluidlimit2g}
q_r^*=z_r^* \Big|_{x_r^*=0}=\beta_r\lambda(1-\lambda_2^*)+\nu \frac{1-\lambda}{R}.
\end{equation}

\cref{zrexpr} implies $z_r^* \leq q_r^*$ and \cref{zrexpr2} leads to $z_r^*\leq \alpha_r\lambda$, $x_r^*=0 \Rightarrow z_r^*=q_r^*$ and $x_r^*>0 \Rightarrow z_r^*=\alpha_r\lambda$, yielding $z_r^*=\min\{q_r^*,\alpha_r\lambda\}$ and thus
 
 \newpage

\begin{equation}\label{exprl1}
\lambda_1^*=\sum_{r=1}^R \min\{q_r^*,\alpha_r\lambda\}.
\end{equation}

%Note that \cref{zrexpr2} leads to the observations $x_r^*>0 \Rightarrow z_r^*=\alpha_r\lambda$ and $x_r^*=0 \Rightarrow z_r^*=q_r^*$. Since $z_r^*$ is maximally $\alpha_r\lambda$, we find $x_r^*>0$ when $q_r^* \geq \alpha_r\lambda$ and $x_r^*=0$ when $q_r^*<\alpha_r\lambda$.
%\cref{fluidlimitexpr4}, with $\beta_r=1/R$ and $\nu=0$ then gives
%\begin{equation}\label{exprl1}
%\lambda_1^*=\sumR \alpha_r\lambda \mathds{1}\{q_r^* \geq \alpha_r \lambda\}+ q_r^* \mathds{1}\{q_r^* < \alpha_r \lambda\}.
%\end{equation}

\emph{Basic JIQ scheme.}
In case $\beta_r=1/R$ and $\nu=0$, 
\cref{exprl1} can be rewritten to
\[
\lambda_2^*=\sum_{r=1}^R\left[\alpha_r\lambda-\frac{\lambda(1-\lambda_2^*)}{R}\right]^+,
\]
and by further calculations, since $\alpha_r$ is decreasing in $r$, to the expression for $\lambda_2^*$ in \cref{maintheorem2}.

In \cref{table3} we compare the fluid-limit approximations
for the mean-waiting time with simulation figures for various
numbers of servers.

\begin{table}[h!!!]
\begin{center}
\caption{Mean waiting time for $\lambda=0.9$ and $\alpha_1=0.8$ or $\alpha_1=0.6$.}
\begin{tabular}{ r c c }
\hline
 & $\lambda=0.9, \alpha_1=0.8$ & $\lambda=0.9, \alpha_1=0.6$\\
N & simulation & simulation\\
\midrule
$10$    & 2.5824 & 2.1234  \\
$20$    & 1.7349 & 1.1386    \\
$50$    & 1.1704 & 0.5001    \\
$100$   & 1.0173 & 0.2981 \\
$200$    & 0.9764 & 0.2207  \\
$500$    & 0.9631 & 0.1983  \\
$1000$ & 0.9599 & 0.1962 \\
\hline
fluid ($N=\infty$) & 0.9643 & 0.1957\\
\hline
\end{tabular}
\label{table3}
\end{center}
\end{table}

%MARKADDED297
\cref{table3} shows that the fluid-limit analysis agrees with the
simulation results, although the number of servers needs to be
larger than in the blocking scenario for extremely high accuracy to
be observed.
Similarly to \cref{table1}, the more symmetric the loads, the better
the performance and the lower the mean waiting time,
which is in line with \cref{maintheorem2}.

\subsection{Enhancements}
We examine the behavior of the system for \cref{alg1,alg2} and show that they can achieve asymptotically zero waiting for any $\lambda<1$ and suitable parameter values as identified in \cref{maintheorem3}. In view of \cref{meanwaitingtime} it suffices to show that $\lambda_2^*=0$ for both enhancements. We first consider \cref{alg1} in which $\alpha_r=\beta_r$ for all $r$ and $\nu=0$. \cref{fluidlimit2g} gives $q_r^* - \alpha_r\lambda = -\beta_r\lambda\lambda_2^*\leq 0$ for all $r$ and $q_r^*=\beta_r\lambda(1-\lambda+\lambda_1^*)$. We obtain
\[
\lambda_1^*=\sum_{r=1}^R q_r^* = \lambda(1-\lambda+\lambda_1^*),
\]
which has a unique solution $\lambda_1^*=\lambda$, so that $\lambda_2^*=0$.\\

Next, we turn to \cref{alg2} where $\nu \geq \frac{\lambda}{1-\lambda}(R \alpha_1-1) $ and $\beta_r=1/R$. If the term in brackets in \cref{zrexpr2} is non-positive, $z_r^*=\alpha_r\lambda$. Otherwise,
\[
\begin{split}
z_r^* &\geq \alpha_r\lambda - \left[ \alpha_r\lambda-\frac{\lambda(1-\lambda_2^*)}{R}-\nu\frac{1-\lambda}{R}\right]\\
%&=\alpha_r\lambda-\alpha_r\lambda + \frac{\lambda - \lambda\lambda_2^*}{R}+\nu\frac{1-\lambda}{R}\\
&\geq\alpha_r\lambda-\alpha_r\lambda + \frac{\lambda - \lambda\lambda_2^*}{R}+\frac{\lambda}{R}(R\alpha_1-1)\geq \alpha_r\lambda - \frac{\lambda\lambda_2^*}{R},
\end{split}
\]
which by \cref{zrexprl1} gives $\lambda-\lambda_2^*=\lambda_1^* \geq \lambda-\lambda\lambda_2^*$, and since $\lambda<1$, we obtain $\lambda_2^*=0$.\\

\begin{figure}[t!!!]
\begin{center}
\includegraphics[width=\scaling\columnwidth]{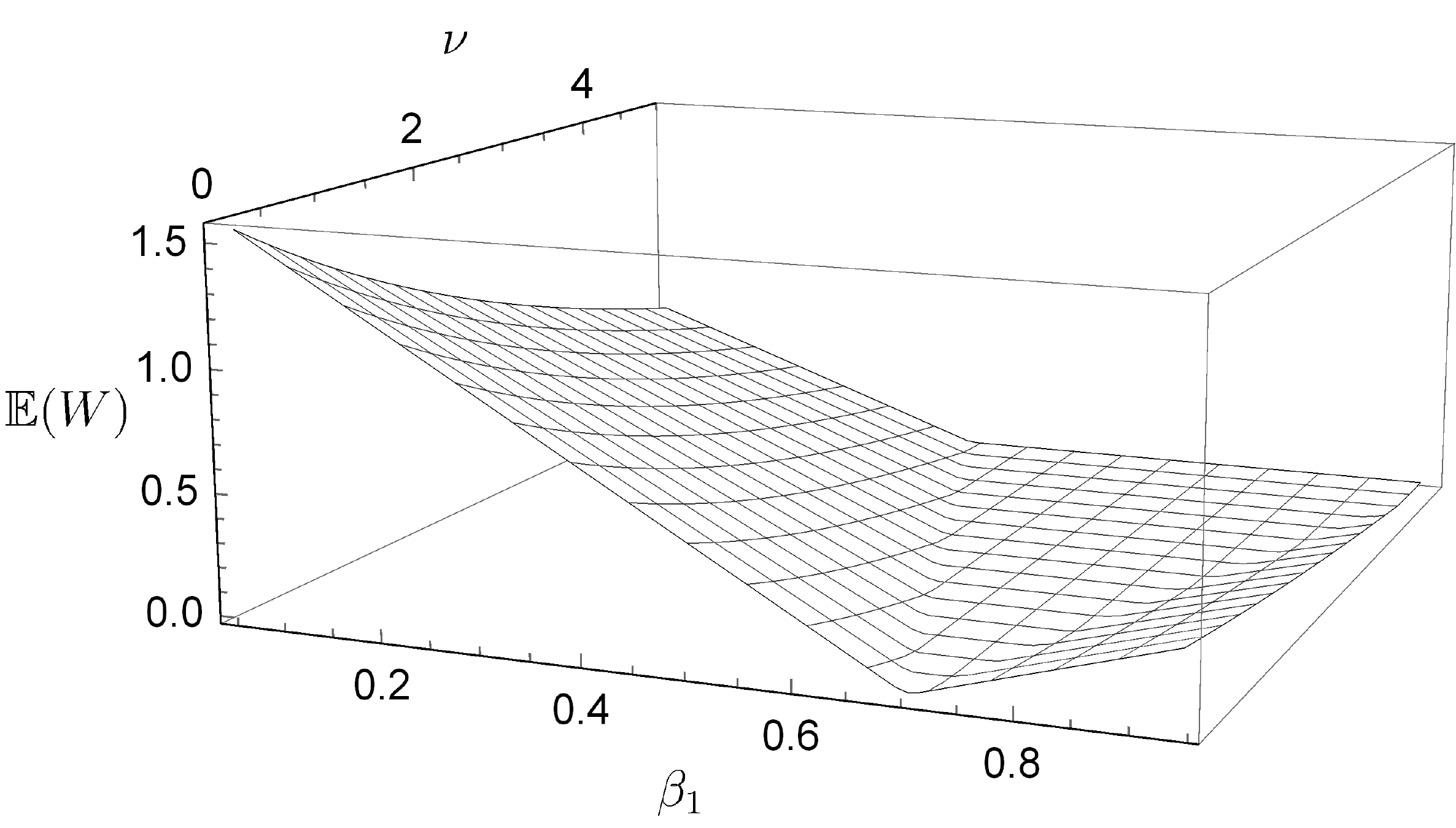}
\caption{Mean waiting time $\mathbb{E}[W]$ in the limit for $R=2$, $\lambda=0.9$ and $\alpha_1=0.7$, for different values of $\beta_1$ and $\nu$.}\label{fwtp4}
\end{center}
\end{figure}

\cref{fwtp4} displays the mean waiting time as $N\to\infty$ for the system with both enhancements. 
%ADDEDMARK287
Similarly to \cref{fwtp3}, we can greatly improve the performance
by tuning~$\beta$ and~$\nu$.
Again $\alpha_1=0.7$, so that $\beta_1=0.7$ is the best choice.
The mean waiting time decreases as $\beta_1$ approaches~$\alpha_1$,
or as the rate~$\nu$ increases.
Exact knowledge of the arrival rates is not required, and a rough
approximation of~$\beta_1$ and a small value of~$\nu$ are
sufficient for the mean waiting time to vanish.

\section{Conclusion}
\label{conc}

We examined the performance of the Join-the-Idle-Queue (JIQ) scheme
in large-scale systems with several possibly heterogeneous dispatchers.
We used product-form representations and fluid limits to show that
the basic JIQ scheme fails to deliver zero blocking and wait for any
asymmetric dispatcher loads, even for arbitrarily low overall load.
Remarkably, it is the least-loaded dispatcher that throttles tokens
and leaves idle servers stranded, thus acting as bottleneck.

In order to counter the performance degradation for asymmetric
dispatcher loads, we introduced two extensions of the basic JIQ
scheme where tokens are either distributed non-uniformly
or occasionally exchanged among the various dispatchers.
We proved that these extensions can achieve zero blocking and wait
in the many-server limit, for any subcritical overall load
and arbitrarily skewed load profiles.
Extensive simulation experiments corroborated these results,
indicating that they apply even in moderately sized systems.

It is worth emphasizing that the proposed enhancements involve no
or constant additional communication overhead per job,
and hence retain the scalability of the basic JIQ scheme.
The algorithms do rely on suitable parameter settings,
and it would be of interest to develop learning techniques for that.

While we allowed the dispatchers to be heterogeneous, we assumed
all the servers to be statistically identical, and the service
requirements to be exponentially distributed. As noted earlier, \cref{maintheorem} in fact holds for non-exponential service requirement distributions as well.
In ongoing work we aim to extend \cref{maintheorem2,maintheorem3} to possibly non-exponential service requirement distributions.

\section*{Acknowledgments}
This work is supported by the NWO Gravitation Networks grant 024.002.003, 
an NWO TOP-GO grant and an ERC Starting Grant.

\appendix

\section{Proof of \cref{maintheorem}}
We now present the analyses of the models required for the proof of \cref{maintheorem}, along the lines sketched in \cref{jacksonnetwork}.
\subsection{Blocking model with two dispatchers
(better system)}
\label{modelA}
Consider the blocking scenario with $N$ servers, $R = 2$ dispatchers
and arrival fractions $\alpha_1$ and $\alpha_2$.
The probability of sending a token to dispatcher~$i$ is~$\beta_i$.
Without loss of generality,
assume $\frac{\alpha_1}{\beta_1} \geq \frac{\alpha_2}{\beta_2}$.
Since this is a closed Jackson network (see \cref{jacksonnetwork}),
we have the stationary distribution
\begin{equation}
\label{eqrp}
\pi(n_0,n_1,n_2)
%=& G^{-1} \frac{(2\beta_1\lambda N)^{n_1}}{(\alpha_1\lambda N)^{n_1}}\frac{(2\beta_2\lambda N)^{n_2}}{(\alpha_2\lambda N)^{n_2}} \frac{(2\lambda N)^{n_0}}{n_0!}\\
= G^{-1} \frac{(\lambda N)^{n_1}}{(\frac{\alpha_1}{2\beta_1}\lambda N)^{n_1}}\frac{(\lambda N)^{n_2}}{(\frac{\alpha_2}{2\beta_2}\lambda N)^{n_2}} \frac{(2\lambda N)^{n_0}}{n_0!},
\end{equation}
with $G$ the normalization constant.
We use this to prove the following proposition.

\begin{myproposition}[Limiting blocking probability for $R=2$]
\label{propA}
\begin{equation*}
\limN B(2, N,\lambda, \hdots) = \max\left\{\beta_1\left(\frac{\alpha_1}{\beta_1}-\frac{\alpha_2}{\beta_2}\right),1-1/\lambda\right\}.
\end{equation*}
\end{myproposition}

Define $\alphas_i=\frac{\alpha_i}{2\beta_i}$. The proof of \cref{propA} starts from an exact expression for the
blocking probability that follows from \cref{eqrp}:
%\begin{equation}\label{derBlocking}
%\begin{split}
%&B(2, N,\lambda, \alpha_1, \alpha_2, \beta_1, \beta_2)\\
%=&(\alphas_1-\alphas_2)
%\frac{\alpha_1 x_2 +\alpha_2 \left(\frac{\alphas_2}{\alphas_1}\right)^N x_1 }{   \alphas_1 x_2 - \alphas_2\left(\frac{\alphas_2}{\alphas_1}\right)^{N} x_1 }
%\end{split}
%\end{equation}
%
%We rewrite this expression as
\begin{equation}
\label{blockingwithZ}
B(2, N, \lambda, \hdots)
=(\alphas_1-\alphas_2)\frac{1+\frac{\alpha_2}{\alpha_1}\frac{\alphas_1}{\alphas_2}Z(N,\lambda,\alphas_1,\alphas_2)}{\frac{\alphas_1}{\alpha_1}-\frac{\alphas_1}{\alpha_1}Z(N,\lambda,\alphas_1,\alphas_2)},
\end{equation}
with
\begin{equation*}
Z(N,\lambda,\alphas_1,\alphas_2)=\left(\frac{\alphas_2}{\alphas_1}\right)^{N+1}\frac{\sumN \frac{\left(2\alphas_1\lambda N \right)^n}{n!}}{\sumN \frac{\left(2\alphas_2\lambda N \right)^n}{n!}}.
\end{equation*}
%with
%\[
%x_i=\sumN \frac{\left(2\alphas_i\lambda N \right)^n}{n!}.
%\]
Through tedious calculations, it may be shown that
\begin{equation*}
\underset{N\rightarrow \infty}{\lim} Z(N,\lambda,\alphas_1,\alphas_2)=
\begin{cases}
0 & 2\alphas_2\lambda\leq 1,\\
\frac{\alphas_2}{\alphas_1}\frac{2\alphas_2\lambda-1}{2\alphas_2\lambda-\frac{\alphas_2}{\alphas_1}} & 2\alphas_2\lambda>1.
\end{cases}
\end{equation*}
Substitution into \cref{blockingwithZ} then yields \cref{propA}.

%%%%%%%%%%
\subsection{Symmetric blocking model (worse system)}
\label{modelB}

Consider the blocking model with $N$ servers, $R$ dispatchers and $\alpha_1=\hdots =\alpha_R=1/R$, assume $\lambda\leq 1$ and use the stationary distribution $\pi$ and blocking probability $B(R,N,\lambda)$ provided in \cref{jacksonnetwork}.

\begin{myproposition}[Recursive formula blocking probability]
\label{recfor}
\begin{equation}\label{reccc}
B(R+1,N,\lambda)=\frac{R}{N+R-\lambda N(1-B(R,N,\lambda))}.
\end{equation}
\end{myproposition}

\begin{myproposition}[Vanishing blocking probability in symmetric systems]
\label{propB}
For $\lambda < 1$, $\underset{N\rightarrow \infty}{\lim} B(R,N,\lambda)=0$.
\end{myproposition}

The proofs again start from an exact expression for the blocking
probability obtained from the stationary distribution:
\begin{equation}\label{blockingprobR}
B(R, N, \lambda)=\frac{\sumN \binom{N+R-2-n}{R-2}\frac{(\lambda N)^n}{n!}}{\sumN \binom{N+R-1-n}{R-1}\frac{(\lambda N)^n}{n!}}.
\end{equation}
\cref{recfor} follows from inserting \cref{blockingprobR} at both
sides of \cref{reccc}. \cref{propB} follows from
\begin{align*}
B(R+1, N,\lambda)&=\frac{R}{R+N(1-\lambda) +\lambda N B(R, N,\lambda)}\nonumber\\
 &\leq  \frac{R}{R+N(1-\lambda)}\goesto 0.
\end{align*}


\begin{thebibliography}{99}

%\bibitem{Asmussen03}
%S. Asmussen (2003).
%{\em Applied Probability and Queues},
%Springer-Verlag, New York.

\bibitem{BB08}
R. Badonnel, M. Burgess (2008).
Dynamic pull-based load balancing for autonomic servers.
{\em Proc.\ IEEE NOMS 2008}, 751--754.
%Salvador, Bahia, Brazil,
%April 7--11, 2008,
%751--754.

\bibitem{BLP10}
M. Bramson, Y. Lu, B. Prabhakar (2010).
Randomized load balancing with general service time distributions.
{\em ACM SIGMETRICS Perf. Eval. Rev.} {\bf 38 (1)},
275--286.

\bibitem{BLP12}
M. Bramson, Y. Lu, B. Prabhakar (2012).
Asymptotic independence of queues under randomized load balancing.
{\em Queueing Systems} {\bf 71 (3)},
247--292.

%\bibitem{BLP13}
%M. Bramson, Y. Lu, B. Prabhakar (2013).
%Decay of tails at equilibrium for FIFO join the shortest queue networks.
%{\em Ann.\ Appl.\ Prob. {\bf 23},
%1841--1878.

%\bibitem{EVW80}
%A. Ephremides, P. Varaiya, J. Walrand (1980).
%A simple dynamic routing problem.
%{\em IEEE Trans.\ Aut. Contr.} {\bf 25},
%690--693.

%\bibitem{Erlang17}
%A.K. Erlang (1917).
%Solution of some problems in the theory of probabilities
%of significance in automatic telephone exchanges.
%In: {\em Elektrotkeknikeren}.

%\bibitem{EG15}
%P. Eschenfeldt, D. Gamarnik(2015).
%Joint the shortest queue with many servers: The heavy traffic asymptotics.
%Preprint, arXiv:1502.00999.

%\bibitem{FS16}
%S. Foss, A.L. Stolyar (2016).
%Large-scale Join-Idle-Queue system with general service times.
%Preprint, arXiv:1605.05968

\bibitem{GTZ16}
D. Gamarnik, J. Tsitsiklis, M. Zubeldia (2016).
Delay, memory and messaging tradeoffs in distributed service systems.
{\em Proc. ACM SIGMETRICS Perf. Eval. Rev.} {\bf 44 (1)}, 1--12.

%\bibitem{GHSW07}
%V. Gupta, M. Harchol-Balter, K. Sigman, W. Whitt (2007).
%Analysis of join-the-shortest-queue routing for web server farms.
%{\em Perf.\ Eval.} {\bf 64 (9--12)},
%1062--1081.

%\bibitem{HW81}
%S. Halfin, W. Whitt (1981).
%Heavy-traffic limits for queues with many exponential servers.
%{\em Oper.\ Res.} {\bf 29},
%567--588.

\bibitem{HK94}
P. Hunt, T. Kurtz (1994).
Large loss networks.
{\em Stoc.\ Proc.\ Appl.} {\bf 53 (2)},
363--378.

%\bibitem{Jagerman74}
%D. L. Jagerman (1974).
%Some properties of the Erlang loss function.
%In: {\em Bell Syst.\ Techn.\ J.} {\bf 53 (3)},
%525--551.

\bibitem{Kelly79}
F. P. Kelly (2011).
Reversibility and stochastic networks.
{\em Cambridge University Press}, New York.

\bibitem{LXKGLG11}
Y. Lu, Q. Xie, G. Kliot, A. Geller, J. Larus, A. Greenberg (2011).
Join-idle-queue: A novel load balancing algorithm for dynamically
scalable web services.
{\em Perf.\ Eval.} {\bf 68 (11)},
1056--1071.

%\bibitem{MS15}
%S. Maguluri, R. Srikant (2015).
%Scheduling jobs with unknown duration in clouds.
%{\em IEEE/ACM Transactions on Networking}, to appear.

%\bibitem{MSY12}
%S. Maguluri, R. Srikant, L. Ying (2012).
%Stochastic models of load balancing and scheduling in cloud
%computing clusters.
%In: {\em Proc.\ IEEE INFOCOM 2012}, Orlando FL, March 27--29, 2012.

%\bibitem{MYS14}
%S. Maguluri, L. Ying, R. Srikant (2013).
%Heavy traffic optimal resource allocation algorithms for cloud
%computing clusters.
%In: {\em Proc.\ ITC 25}, Shanghai, China, September 10--12, 2013.

\bibitem{Mitzenmacher01}
M. Mitzenmacher (2001).
The power of two choices in randomized load balancing.
{\em IEEE Trans.\ Par.\ Distr.\ Syst.} {\bf 12 (10)},
1094--1104.

\bibitem{Mitzenmacher16}
M. Mitzenmacher (2016).
Analyzing distributed Join-Idle Queue: A fluid limit approach.
Preprint, arXiv:1606.01833.

%\bibitem{MBLW16a}
%D. Mukherjee, S.C. Borst, J.S.H. van Leeuwaarden, P.A. Whiting (2016).
%Universality of load balancing schemes on diffusion scale.
%{\em J. Appl.\ Prob.}, to appear.

\bibitem{MBLW16c}
D. Mukherjee, S.C. Borst, J.S.H. van Leeuwaarden, P.A. Whiting (2016).
Asymptotic optimality of power-of-$d$ load balancing in large-scale systems.
Preprint, arXiv:1612.00722.

\bibitem{MBLW16d}
D. Mukherjee, S.C. Borst, J.S.H. van Leeuwaarden, P.A. Whiting (2016).
Universality of power-of-$d$ load balancing in many-server systems.
Preprint, arXiv:1612.00723.

\bibitem{MBLW16b}
D. Mukherjee, S.C. Borst, J.S.H. van Leeuwaarden, P.A. Whiting (2016).
Universality of power-of-$d$ load balancing schemes.
{\em ACM SIGMETRICS Perf. Eval. Rev.}  {\bf 44 (2)}, 36--38.

\bibitem{MKM15}
A. Mukhopadhyay, A. Karthik, R.R. Mazumdar (2015).
Randomized assignment of jobs to servers in heterogeneous clusters
of shared servers for low delay.
{\em Stoc. Syst.} {\bf 6 (1)}, 90--131.

\bibitem{MKMG15}
A. Mukhopadhyay, A. Karthik, R.R. Mazumdar, F. Guillemin (2016).
Mean field and propagation of chaos in multi-class heterogeneous
loss models.
{\em Perf.\ Eval.} {\bf 91},
117--131.

%\bibitem{MM14a}
%A. Mukhopadhyay, R.R. Mazumdar (2014).
%Randomized routing schemes for large processor sharing systems
%with multiple service rates.
%In: {\em Proc.\ ACM SIGMETRICS 2014}.
%Austin TX, June 17--19.

%\bibitem{MM14b}
%A. Mukhopadhyay, R.R. Mazumdar (2014).
%Rate-based randomized routing in large heterogeneous processor
%sharing systems.
%In: {\em Proc.\ ITC 26}, Karlskrona, Sweden, September 9--11, 2014.

%\bibitem{MMG15}
%A. Mukhopadhyay, R.R. Mazumdar, F. Guillemin (2015).
%The power of randomized routing in heterogeneous loss systems.
%In: {\em Proc.\ ITC-27}.

%\bibitem{ananta}
%P. Patel, D. Bansal, L. Yuan, A. Murthy, A. Greenberg, D.A. Maltz,
%R. Kern, H. Kumar, M. Zikos, H. Wu, C. Kim, N. Karri (2013).
%Ananta: Cloud-scale load balancing.
%{\em ACM SIGCOMM Comp.\ Commun.\ Rev.} {\bf 43 (4)},
%207--218.

\bibitem{Stolyar15a}
A.L. Stolyar (2015).
Pull-based load distribution in large-scale heterogeneous service systems.
{\em Queueing Systems} {\bf 80 (4)},
341--361.

\bibitem{Stolyar15b}
A.L. Stolyar (2015).
Pull-based load distribution among heterogeneous parallel servers:
the case of multiple routers.
{\em Queueing Systems}, to appear.

\bibitem{VDK96}
N. Vvedenskaya, R. Dobrushin, F. Karpelevich (1996).
Queueing system with selection of the shortest of two queues:
An asymptotic approach.
{\em Prob.\ Inf.\ Trans.} {\bf 32 (1)},
20--34.

%\bibitem{Weber78}
%R.R. Weber (1978).
%On the optimal assignment of customers to parallel queues.
%{\em J.\ Appl.\ Prob.\} {\bf 15},
%406--413.

%\bibitem{Whitt86}
%W. Whitt (1986).
%Deciding which queue to join: Some counterexamples.
%{\em Oper.\ Res.} {\bf 34 (1)},
%55--62.

%\bibitem{Winston77}
%W. Winston (1977).
%Optimality of the shortest line discipline.
%{\em J. Appl.\ Prob.} {\bf 14}, 181--189.

\bibitem{XDLS15}
Q. Xie, X. Dong, Y. Lu, R. Srikant (2015).
Power of $d$ choices for large-scale bin packing: A loss model.
{\em ACM SIGMETRICS Perf. Eval. Rev.} {\bf 43 (1)}, 321--334.
%Portland OR, June 16--18, 2015.

\bibitem{YSK15}
L. Ying, R. Srikant, X. Kang (2015).
The power of slightly more than one sample in randomized load balancing.
{\em Proc.\ IEEE INFOCOM 2015}, 1131--1139.
%Hong-Kong, April 28--30, 2015.


\end{thebibliography}
\end{document}